\documentclass[a4paper]{article}
\usepackage[latin1]{inputenc}
\usepackage{graphicx}
\usepackage{geometry}
\usepackage{amssymb, amsmath, amsfonts, amsthm}
\usepackage{bbm}
\usepackage{MnSymbol}
\usepackage{mathrsfs}

\usepackage{url}
\usepackage{hyperref}
\usepackage{color}
\usepackage{euscript}

\newcommand{\leaveout}[1]{\centerline{
\framebox{ Something has been left out here}}}

\renewcommand{\d}{\ensuremath{\mathrm{d}}}
 
\newcommand{\db}{\ensuremath{\bar{\mathrm{d}}}}

\newcommand{\g}{\ensuremath{\mathfrak{g}}}
\newcommand{\R}{\ensuremath{\mathbb{R}}}

\newcommand{\dg}{\ensuremath{\overline{\nabla}}}

\newtheorem{theorem}{Theorem}

\newtheorem{example}[theorem]{Example}

\newcommand{\grad}{\ensuremath{\mathrm{grad}}}
\newcommand{\ret}{\phi}

\newcommand{\intder}{\ensuremath{\righthalfcup}}

\title{\bf Preserving first integrals with symmetric Lie group methods  }
\author{E. Celledoni and B. Owren}

\begin{document}
\maketitle
\begin{center}
{\it In honour of Arieh Iserles}
\end{center}

\begin{abstract}
The discrete gradient approach is generalized to yield integral preserving methods for differential equations in Lie groups. 
 \end{abstract}
 
\section{Introduction}
Our point of departure is the system of differential equations
\begin{equation} \label{eq:diffeq}
      \dot{x} = F(x) = f(x)\cdot x = x\cdot \tilde{f}(x),
\end{equation}
where the unknown $x=x(t)$ is a curve on some Lie group $x(t)\in G$ and the dot over $x$ signifies differentiation with respect to $t$.
The map $f: G\rightarrow\g$ where $\g\equiv T_eG$ is the Lie algebra corresponding to $G$, and the dot should be interpreted as the derivative of right (resp. left) multiplication, e.g.
$f(x)\cdot x := T_eR_x f(x)$. Such equations occur for instance in mechanical systems where the Lie group could be the special orthogonal group $SO(3)$ or the special  Euclidean group $SE(3)$. We shall consider in particular the case where the system
\eqref{eq:diffeq} possesses one or more first integrals. Here we define a first integral to be any function $H: G\rightarrow \R$ which is invariant on solutions 
$$
       \frac{\mathrm{d}}{\mathrm{d}t}H(x(t)) = \langle\d H, F\rangle = 0.
$$
Any differential equation \eqref{eq:diffeq} on a Lie group having $H$ as first integral can be formulated via a bivector (dual two-form)
$\omega$ and the differential of $H$ as
\begin{equation} \label{eq:diffeqbivector}
    \dot{x} = F(x) = \omega(\d H,\cdot) = \d H \intder \omega
\end{equation}
A Riemannian metric on $G$ defines an inner product $(\cdot,\cdot)_x$ on every tangent space $T_xG$ which is also varying smoothly with $x$. 
With such a metric one may define the Riemannian gradient vector field as the unique vector field $\grad H$ satisfying
$\langle\d H, v\rangle_x=(\grad H|_x,v)_x$ for every $v\in T_xG$.
 An example of a bivector $\omega$ to be used in \eqref{eq:diffeqbivector} is then provided by means of the wedge product
$$
      \omega = \frac{\grad H\wedge F}{\|\grad H\|^2}.
$$
Note that $\omega$ is not uniquely defined by $F$ and $H$. For a choice $x=(x_1,\ldots,x_d)$ of local coordinates on the group, we may write
\begin{equation} \label{eq:bivectorcoords}
    \omega|_x = \sum_{1\leq i<j\leq d} \omega_{i,j}(x)\,\partial x_{i} \wedge \partial x_j,\qquad
    \d H|_x = \sum_{k=1}^d \frac{\partial H}{\partial x_k}\,\d x_k
\end{equation}
In this way, we find the coordinate version of \eqref{eq:diffeqbivector}
\begin{equation} \label{eq:coordvs}
\dot{x} = S_\omega \nabla H
\end{equation}
for the skew-symmetric $d\times d$-matrix $S_\omega=\boldsymbol\omega^T-\boldsymbol\omega$ and $\boldsymbol\omega$ is the strictly upper triangular matrix with
entries $\omega_{ij}(x),\ 1\leq i<j\leq d$. The formulation \eqref{eq:coordvs} has been the starting point for energy preserving integrators devised by Gonzalez \cite{gonzalez96tia}. In fact, McLachlan et al. \cite{mclachlan99giu} showed that under relatively general circumstances, any vector field $F$ on $\mathbb{R}^d$ with a first integral $H$ can be written in the form
\begin{equation} \label{eq:mclachlan}
      F(x) = S(x)\nabla H
\end{equation}
for some skew-symmetric matrix $S(x)$.

Note that the bivector is not required to be nondegenerate. For Hamiltonian systems, the Hamiltonian itself is a first integral and an
 accompanying bivector can be inferred from the symplectic two-form $\Omega$ by inversion.
In coordinates one can represent the symplectic two-form by a skew-symmetric matrix $S_\Omega$ in a similar way
as in \eqref{eq:bivectorcoords} and \eqref{eq:coordvs} with respect to the basis $\{\d x_i \wedge \d x_j,\ i<j\}$. By definition, this matrix will be
invertible, and its inverse is precisely $S_\omega=S_\Omega^{-1}$.

The formulation \eqref{eq:diffeqbivector} is easily generalised to the case with $k$ independent first integrals $H_1,\ldots, H_k$. We may now replace the
bivector with a $k+1$-vector $\omega$ and write \eqref{eq:diffeq} as
$$
    \dot{x} = F(x)=\omega(\d H_1,\ldots, \d H_k,\;\cdot\;)
$$
Also in this case we can find, by means of a Riemannian structure, an example of a feasible $k+1$-vector
$$
    \omega = \frac{\omega_0\wedge F}{\omega_0(\d H_1,\ldots,\d H_k)}\quad\mbox{where}\quad
    \omega_0 = \grad H_1\wedge\cdots\wedge\grad  H_k
$$
Earlier work on energy-preserving methods on Lie-groups was presented in \cite{lewis94caf}, see also the references therein.
In this paper we shall generalise the notion of discrete gradient methods to the situation where the phase space is a Lie group.
The paper extends results presented in \cite{celledoni12ait} by also considering high order methods, a larger class of manifolds and further examples. 

\section{Discrete differentials in Lie groups}
\subsection{A review of the situation in Euclidean space}
The idea of discrete gradient methods, is to consider some approximation to the exact gradient in \eqref{eq:mclachlan}
$\dg H: \mathbb{R}^d\times\mathbb{R}^d\rightarrow\mathbb{R}^d$ satisfying the following two conditions
\begin{align}
H(v)-H(u)&=\dg H(u,v)^T (v-u),&\forall u,v\in\mathbb{R}^d\label{eq:dchainrule} \\
\dg H(u,u)&=\nabla H(u),&\forall u\in\mathbb{R}^d \label{eq:dconsistent}
\end{align}
Many such discrete gradients have been proposed in the literature, and we give here a few examples.
The averaged vector field discrete gradient is discussed for instance in \cite{mclachlan99giu} and more recently in the PDE setting 
\cite{dahlby11agf, celledoni12per}
\begin{equation} \label{eq:cavf}
       \dg H(u,v) = \int_0^1 \nabla H((1-\xi)u+\xi v)\,\d\xi.
\end{equation}
Another discrete gradient which was proposed in \cite{gonzalez96tia} is the midpoint gradient
\begin{equation}\label{eq:ogmp}
       \dg H(u,v) = \nabla H\left(\frac{u+v}{2}\right) + \frac{H(v)-H(u)-\nabla H\left(\frac{u+v}{2}\right)^T(v-u)}{\|v-u\|^2}(v-u)
\end{equation}
An integral preserving integrator for \eqref{eq:mclachlan} is readily given as
\begin{equation} \label{eq:cmethod}
     \frac{x^{n+1}-x^n}{h} = \bar{S}(x^n,x^{n+1})\dg H(x^n,x^{n+1})
\end{equation}
where $\bar{S}$ is a skew-symmetric matrix which approximates $S$. It is required to satisfy the consistency condition $\bar{S}(u,u)=S(u)$.
An easy calculation, using \eqref{eq:dchainrule} and \eqref{eq:cmethod} now shows that
$$
   H(x^{n+1})-H(x^n) = \dg H(x^n,x^{n+1})(x^{n+1}-x^n) =h \dg H(x^n,x^{n+1}) \bar{S}(x^n,x^{n+1})\dg H(x^n,x^{n+1})=0,
$$
the last identity follows since $\bar{S}$ is skew-symmetric. The skew-symmetric matrix $\bar{S}$ used in the integrator is not unique, and the freedom can be used to improve the approximation \eqref{eq:cmethod} in various ways, for instance to increase its order of convergence.
\subsection{The Lie group setting}
For Lie groups, none of the definitions \eqref{eq:dchainrule} or \eqref{eq:cmethod} can be used since they both involve vector space operations, not generally defined on Lie groups. In a coordinate free setting, we also find it more convenient to replace the 
gradient and its discrete counterpart by dual quantities, the differential as indicated in \eqref{eq:diffeqbivector}. Below, we also use the notion of a \emph{discrete differential} rather than a discrete gradient.
As is the tradition for Lie group integrators \cite{iserles00lgm,christiansen11tis}, approximations are introduced through some finite dimensional action and a trivialisation principle is applied. By this, we mean that tangent vectors at some $x\in G$, i.e. $v_x\in T_xG$ can be represented via either left or right translation of a vector $\xi\in T_eG\cong \g$, where $\g$ is the Lie algebra of $G$. In the present paper we will just for convenience choose right translation.
Defining the right multiplication operator $R_x:G\rightarrow G$, where $R_x y = y\cdot x$, we shall use the notation
$$
     v_x = T_eR_x\xi =: R_{x*}\xi
$$
and similarly,  any vector $p_x\in T_x^*G$ is identified by some $\mu\in\g^*$ through
$$
     \langle p_x, R_{x*}\xi\rangle = \langle R_{x}^*p_x, \xi\rangle = \langle\mu,\xi\rangle\quad\mbox{i.e.}\quad \mu=R_x^*p_x,
$$
where $\langle\cdot,\cdot\rangle$ is some duality pairing between $T_x^*G$ and $T_xG$, as well as between $\g^*$ and $\g$.
For example, if the Lie group $G$ and its Lie algebra $\g$ both are realized as $m\times m$-matrices, the dual elements can also be represented as matrices and the duality pairing could be given as $\langle p, v\rangle=\mbox{trace}(p^Tv)$. In this case one simply has
$v_x=R_{x*}\xi = \xi\cdot x$, and $\mu=R_x^*p_x =p_x\cdot x^T$. We also note that in Euclidean space
($G=\mathbb{R}^d$ and group operation is $+$), left and right translation are both realised by the identity map, e.g.
$R_{x*}\xi = \xi$.

We shall introduce the \emph{trivialised discrete differential} of a function $H$ as a map $\db H: G\times G\rightarrow\g^*$ satisfying the following identities generalised from \eqref{eq:dchainrule} and \eqref{eq:dconsistent}
\begin{align}
H(v)-H(u) &= \langle \db H(u,v), \log(v\cdot u^{-1})\rangle \label{eq:chainrule} \\
\db H(x,x) &= R_x^*\d H_x \label{eq:consistency}
\end{align}
By identifying vectors and co-vectors in Euclidean space through the standard inner product, and noting that the logarithmic map
$\log: G\rightarrow\g$ in this setting is simply the identity map, thus $\log (v\cdot u^{-1})=v-u$, we recover the standard discrete gradient conditions 
 \eqref{eq:dchainrule} and \eqref{eq:dconsistent} when Euclidean space is chosen as our Lie group.
 We also need to introduce a trivialised approximation to the bivector $\omega$ in \eqref{eq:diffeqbivector}. For this purpose we define, for any pair of points $(u,v)\in G\times G$, an exterior 2-form on the linear space $\g^*$ which we denote by $\bar{\omega}(u,v)$ thus
 $\bar{\omega}: G\times G\rightarrow \Lambda^2(\g^*)$. We impose the consistency condition 
  \begin{equation*}
     \bar{\omega}(x,x)(R_x^*\alpha, R_x^*\beta) = \omega_x(\alpha,\beta),\qquad \forall x\in G,\quad\forall \alpha,\beta\in T_x^*G
 \end{equation*}
 Of course, in practice, $\bar{\omega}$ need only be defined in some suitable neighborhood of the diagonal subset $\{(x,x),\ x\in G\}$.
 Introducing coordinates, the form $\bar{\omega}$ plays a similar role as the skew-symmetric matrix
 $\bar{S}(x^n,x^{n+1})$ in \eqref{eq:cmethod}. We may further define our numerical method as follows
 \begin{equation} \label{eq:methdef}
    x^{n+1} = \exp(h F(x^n,x^{n+1}))\cdot x^n,\qquad F(x^n,x^{n+1}) = \db H(x^n,x^{n+1}) \righthalfcup \bar{\omega}(x^n,x^{n+1})
 \end{equation}
 For the reader who is unfamiliar with the notation, the definition of $F(x^n,x^{n+1})\in\g$ using coordinates as in \eqref{eq:coordvs}
 would be $F(x^n,x^{n+1}) = \bar{S}_\omega(x^n,x^{n+1}) \db H(x^n,x^{n+1})$ where $\bar{S}_\omega(x^n,x^{n+1})$ is a 
 skew-symmetric matrix approximating $S_\omega(x^n)$ and where we have also expressed $\db H(x^n,x^{n+1})$ in coordinates.
 
 From the defining relation \eqref{eq:chainrule}, it follows immediately that the method preserves $H$ since
\begin{equation*}
 \begin{split}
     H(x^{n+1})-H(x^n)&=\langle \db H(x^n,x^{n+1}), \log(x^{n+1}(x^n)^{-1})\rangle
     =h\langle \db H(x^n,x^{n+1}), F(x^n,x^{n+1})\rangle \\
     &= \bar{\omega}(x^n,x^{n+1})(\db H(x^n,x^{n+1}),\db H(x^n,x^{n+1})\rangle=0
\end{split}
\end{equation*}

\subsection{Examples of trivialised discrete differentials}
An example of a trivialised discrete differential, generalising \eqref{eq:cavf} is
\begin{equation} \label{eq:avf}
\db H(u,v) = \int_0^1 R_{\ell(\xi)}^*\d H_{\ell(\xi)}\d\xi,\qquad \ell(\xi)=\exp(\xi\log(v\cdot u^{-1}))\cdot u
\end{equation}
Here, we have introduced a straight line in the Lie algebra between the points $0$ and $a=\log(v\cdot u^{-1})$.
By applying $\exp$ to each point on the curve and multiplying by $u$, we obtain a curve on the Lie group between the points $u$ and $v$, this is $\ell(\xi)$. Finally the (trivialisation of the)  differential $\d H$ is averaged along this curve to obtain the AVF type of trivialised discrete differential.
Considering $\ell(\xi)$ with $u$ and $v$ interchanged yields $\ell(1-\xi)$, and from this it easily follows that
$\db H(u,v)=\db H(v,u)$.

The Gonzalez midpoint gradient can be generalised as well, for instance by introducing an inner product on the Lie algebra, we denote it
$(\cdot,\cdot)$. We apply ``index lowering" to any element $\eta\in\g$ by defining $\eta^\flat\in\g^*$ to be  the unique element satisfying
$\langle\eta^\flat, \zeta\rangle=(\eta,\zeta)$ for all $\zeta\in\g$.  We can then introduce a generalisation of \eqref{eq:ogmp} as
\begin{equation} \label{eq:tddgmp}
\db H(u,v) = R_c^*\d H|_c + \frac{H(v)-H(u)-\langle R_c^* \d H|_c,\eta\rangle}{(\eta,\eta)}\,\eta^\flat,\quad
\eta = \log(v\cdot u^{-1}),
\end{equation}
where $c\in G$, is some point typically near $u$ and $v$. One may for instance choose $c=\exp(\eta/2)\cdot u$, which implies symmetry, i.e. $\db H(u,v)=\db H(v,u)$.

We have now presented two examples of trivialised discrete differentials which are both symmetric in the two arguments $u$ and $v$.
One observes from the definition \eqref{eq:methdef} that the integrator itself is symmetric if $\db H(u,v)=\db H(v,u)$ and
$\bar{\omega}(u,v)=\bar{\omega}(v,u)$ for all pairs $(u,v)$.

\section{More general manifolds}
More interesting examples of mechanical systems can be found for instance  in the larger class of homogeneous manifolds.
There are various ways to devise discrete gradient methods in this setting, or for even more general classes of manifolds.
From now on, we assume that $M$ is a smooth manifold for which  a retraction map is available, retracting the tangent bundle $TM$ into $M$.
This is a very basic and straightforward approach, certainly there are other ways to devise integral preserving numerical schemes.
A retraction is a map
$$
     \ret: TM\rightarrow M.
$$
Denote by $\ret_p$ the restriction of $\ret$ to $T_pM$ and let $0_p$ be the zero-vector in $T_pM$.
Following \cite{adler02nmo}, we impose the following conditions on $\ret$
\begin{enumerate}
\item $\ret_p$ is smooth and defined in an open ball $B_{r_p}(0_p)\subset T_pM$ of radius $r_p$ about $0_p$.
\item $\ret_p(v)=x$ if and only if $v=0_p$.
\item $T_{0_p}\ret_p=\mathrm{Id}_{T_pM}$. 
\end{enumerate}
This implies in particular that $\ret_p$ is a diffeomorphism from some neighbourhood $\mathscr{U}$ of $0_p$ to 
its image $\mathscr{W}=\ret_p(\mathscr{U})\subset M$. In what follows, we shall always assume that the step size used in the integration is sufficiently small such that both the initial and terminal point of the step are contained in such a set $\mathscr{W}$.
Furthermore, we shall assume that
 there is a given map $c$ defined on some open subset of $M\times M$ containing all
diagonal points $(p,p)$, for which $c(p,q)\in M$. Typically $c(p,q)$ will be either $p$ or $q$ or some kind of centre point between $p$ and $q$ to be defined later. We shall always require that $c(p,p)=p$ for any $p\in M$.

We assume as before the existence of a bivector $\omega$ on $M$ and a first integral $H$ such that the differential equation can be written in the form \eqref{eq:diffeqbivector}. We introduce, for any pair of points $p$ and $q$ on $M$, an approximate bivector  $\bar{\omega}(p,q)$ such that
$$
      \bar{\omega}(p,p)(v,w)=\omega|_p(v,w),\qquad\forall v, w\in T_pM.
$$
The \emph{discrete differential} of a function $H$ can now be defined for any pair of points $(p,q)\in M\times M$ as
a covector $\db H(p,q)\in T_{c(p,q)}^*M$ satisfying the relations
\begin{align*}
   H(q)-H(p) &=  \langle \db H(p,q), \ret_c^{-1}(q)-\ret_c^{-1}(p)\rangle\\
    \db H(p,p) &= \mathrm{d}H|_p,\quad \mbox{for every}\ p\in M. 
\end{align*}
where $c = c(p,q)$ is the map referred to above. We define the integrator  as
\begin{equation} \label{eq:manimeth}
     x^{n+1} = \ret_c(W(x^n,x^{n+1})),\quad W(x^n,x^{n+1})=\ret_c^{-1}(x^n)+h\, \db H(x^n,x^{n+1}) \intder \bar{\omega}(x^n,x^{n+1})
\end{equation}
It follows immediately that the method is symmetric if the following three conditions are satisfied:
\begin{enumerate}
\item The map $c$ is symmetric, i.e. $c(p,q)=c(q,p)$ for all $p$ and $q$.  \label{csym}
\item The discrete differential is symmetric in the sense that $\db H(p,q)=\db H(q,p)$. \label{dbsym}
\item The bivector $\bar{\omega}$ is symmetric in $p$ and $q$: $\bar{\omega}(p,q)=\bar\omega(q,p)$. \label{omsym}
\end{enumerate}
The condition \ref{csym}) can be achieved by solving the equation
\begin{equation}\label{eq:symccond}
     \ret_c^{-1}(p) + \ret_c^{-1}(q) = 0,
\end{equation}
with respect to $c$.

We can now write down a version of the AVF type discrete differential. Let $\gamma_\xi=(1-\xi)v+\xi w$ where
$p=\ret_c(v)$, $q=\ret_c(w)$. Then
\begin{equation} \label{eq:avfgenman}
\db H(p,q) = \int_0^1\ret_c^* \left.\d H\right|_{\ret_c(\gamma_\xi)} \, \d\xi
\end{equation}
Similarly, assuming that $M$ is Riemannian, we can define the following counterpart to the Gonzalez midpoint discrete gradient
\begin{equation} \label{eq:gonzgenman}
\db H(p,q) = \d H|_c + \frac{H(q)-H(p)-\langle \d H|_c,\eta\rangle}{(\eta,\eta)_c}\,\eta^\flat,\quad
 \eta = \ret_c^{-1}(q)-\ret_c^{-1}(p) \in T_cM.
\end{equation}
where we may require that $c(p,q)$ satisfies \eqref{eq:symccond} for the method to be symmetric.

\begin{example}
We consider the sphere $M=S^{n-1}$ where we represent its points as vectors in $\mathbb{R}^n$ of unit length,
$\|p\|_2=1$. The tangent space at $p$ is then identified with the set of vectors in $\mathbb{R}^n$ orthogonal to $p$ with respect to the Euclidean inner product $(\cdot,\cdot)$.
 A retraction is
\begin{equation} \label{eq:rets2}
     \ret_p(v_p) = \frac{p+v_p}{\|p+v_p\|},
\end{equation}
its inverse is defined in the cone $\{q: (p,q)>0\}$ where
$$
     \ret_p^{-1}(q) = \frac{q}{( p, q)} - p
$$
A symmetric map $c(p,q)$ satisfying \eqref{eq:symccond} is simply
\begin{equation} \label{eq:centerpoint}
      c(p,q) = \frac{p+q}{\|p+q\|_2},
\end{equation}
the geodesic midpoint between $p$ and $q$ in terms of the standard Riemannian metric on $S^{n-1}$.
We compute the tangent map of the retraction to be
$$
      T_u\ret_c = \frac{1}{\|c+u\|_2}\left(I-\frac{(c+u)\otimes(c+u)}{\|c+u\|_2^2}\right)
$$
\begin{figure}
\centering
\includegraphics[width=0.5\linewidth]{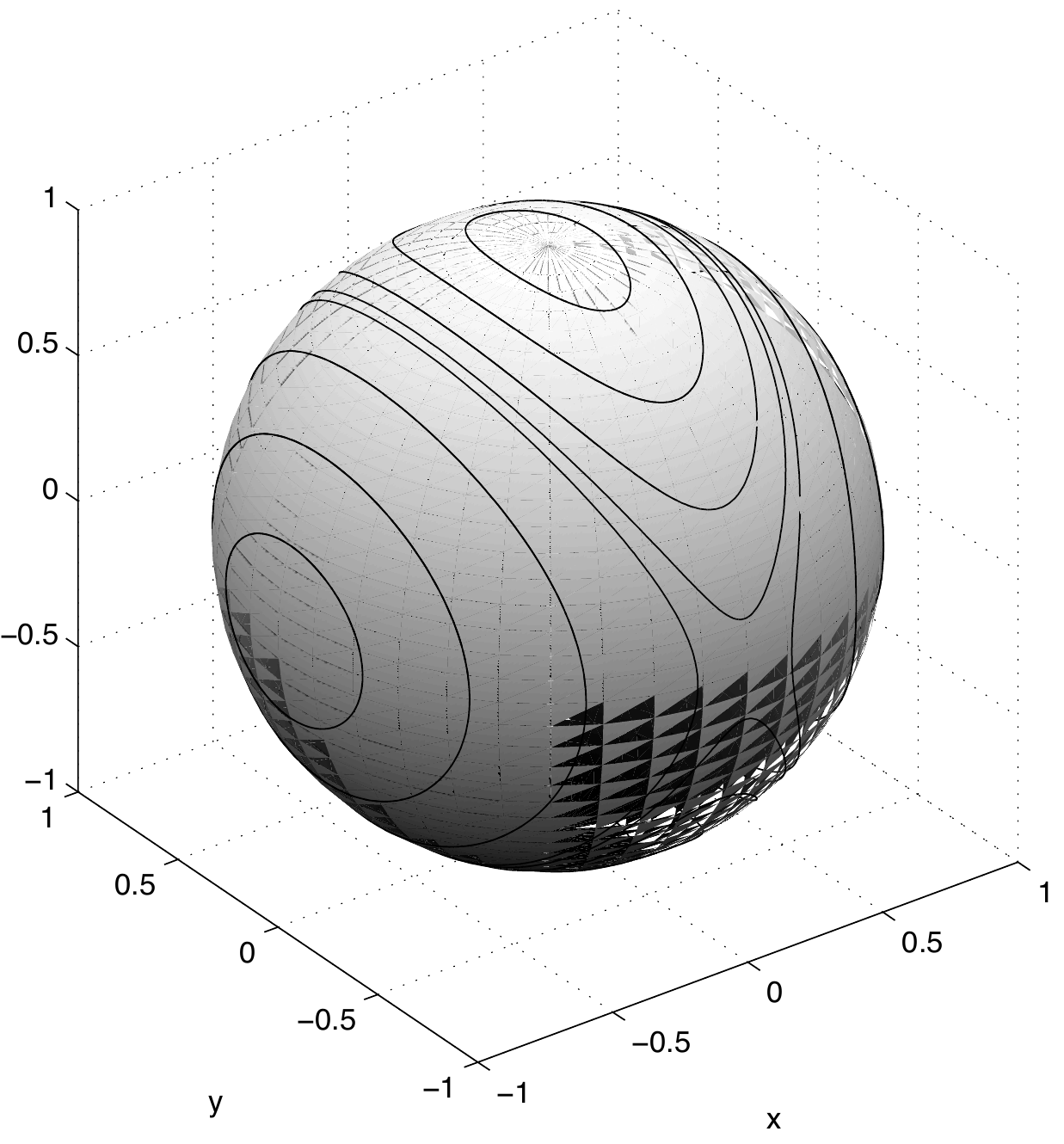}
\caption{The free rigid body equation and its angular momentum in body coordinates. The trajectories are curves of constant energy on the sphere, computed by a method presented here. Moments of inertia used are $\mathbb{I}=\mbox{diag}(1,2,3)$.\label{fig:sphere}}
\end{figure}
As a toy problem, let us consider a mechanical system on $S^2$. Since the angular momentum in body coordinates for the free rigid body is of constant length, we may assume $(p,p)=1$ for all $p$ and we can model the problem as a dynamical system on the sphere. But in addition to this, the energy of the body i preserved,
i.e.
$$
      H(p) = \frac12(p, \mathbb{I}^{-1}p) = \frac12\left(\frac{p_1^2}{\mathbb{I}_1}+\frac{p_2^2}{\mathbb{I}_2}+\frac{p_3^2}{\mathbb{I}_3}\right)
$$
which we may take as the first integral to be preserved. Here the inertia tensor is $\mathbb{I}=\mbox{diag}(\mathbb{I}_1, \mathbb{I}_2, \mathbb{I}_3)$. The system of differential equations can be written as follows
\begin{align*}
\frac{\d p}{\d t} &=  \left.(\d H \intder \omega)\right|_p = p \times \mathbb{I}^{-1} p \\
\omega|_p(\alpha,\beta) &= (p,\alpha\times\beta)
\end{align*}
where the righthand side in both equations refer to the representation in $\mathbb{R}^3$.
A symmetric consistent approximation to $\omega$ would be
$$
  \bar{\omega}(p,q)(\alpha,\beta) =(\frac{p+q}{2}, \alpha\times\beta)
$$

We write $\ell_\xi=c+\gamma_\xi$ with the notation in \eqref{eq:avfgenman}, this is a linear function of the scalar argument $\xi$.
and thus, $\ret_c(\gamma_\xi)=\ell_{\xi}/\|\ell_\xi\|$ from \eqref{eq:rets2}. We therefore derive for the AVF discrete gradient
 $$
\db H(p,q) = \int_0^1 \frac{1}{\|\ell_\xi\|}\left(\mathbb{I}^{-1}\ret_c(\gamma_\xi)  - (\ret_c(\gamma_\xi),\mathbb{I}^{-1}\ret_c(\gamma_\xi))
\ret_c(\gamma_\xi)\right)\;
\d\xi
$$
This integral is somewhat complicated to solve analytically. Instead, we may consider the discrete gradient \eqref{eq:gonzgenman}
where we take as Riemannian metric the standard Euclidean inner product restricted to the tangent bundle of $S^2$.
We obtain the following version of the discrete differential in the chosen representation
$$
\db H(p,q)=  \frac{1}{\|m\|}\left(\mathbb{I}^{-1}m + \frac{\|m\|^2-1}{\|q-p\|^2}(H(q)-H(p))(q-p)\right),\quad m=\frac{p+q}{2}.
$$
The corresponding method is symmetric, thus of second order, and in Figure~\ref{fig:sphere} we used this method to draw the trajectories of the free rigid body problem.

\end{example}

\section{Methods of higher order}
\label{sec:4}


One viable way to obtain higher order variants of the proposed methods is by following the collocation strategy proposed by Hairer in \cite{hairer10epc} and Cohen and Hairer in \cite{cohen11lep}, see also \cite{brugnano10hbv}, and generalise it to the Lie group setting.

Let $c_1,\dots, c_s$ be distinct real numbers such that $0\le c_i\le 1$ and $b_i\ne 0$ for all $i$. We consider  the polynomial $\sigma(\xi h)$ of degree $s$ such that 
\begin{eqnarray}
\label{eq:collmethod1}
\sigma(0)&=&0,\\
\label{eq:collmethod2}
\left.\frac{d}{d \xi} {\sigma}(\xi h)\right|_{\xi=c_j}&=&\mathrm{dexp}_{\sigma_j}^{-1}\left( \bar{d}H_j   \righthalfcup\bar{\omega}_j \right),\quad \sigma_j:=\sigma(c_jh), 
\end{eqnarray}
where
\begin{equation}
\label{eq:DGho}
\bar{d}H_j :=\int_0^1\frac{\ell_j(\xi)}{b_j}\left( \mathrm{dexp}_{\sigma_j}^{-1} \right)^*\mathrm{dexp}_{\sigma(\xi h)}^*\left( R_{\exp(\sigma(\xi h))x_0}^{*}\,dH_{\exp(\sigma(\xi h))x_0}\right)\,d\xi,
\end{equation}
and  $\bar{\omega}_j:\mathfrak{g}^*\times \mathfrak{g}^*\rightarrow \mathbf{R}$ are exterior 2-forms satisfying
\begin{equation}
\label{eq:omegaj}
\bar{\omega}_j(a, b):=\omega_{X_j}(R_{X_j^{-1}}^*\,a,R_{X_j^{-1}}^*\,b),\quad a,b\in\mathfrak{g}^*,
\end{equation}
with $X_j:=R_{x_0} \exp(\sigma_j)$. The numerical solution after one step is
$x_1=R_{x_0}\exp(\sigma(h))$. 

The collocation polynomial is obtained by integrating
\begin{equation}
\label{eq:collpoly}
\frac{d}{d \xi} {\sigma}(\xi h)=\sum_{j=1}^s\ell_j(\xi)\mathrm{dexp}_{\sigma_j}^{-1}\left( \bar{d}H_j   \righthalfcup\bar{\omega}_j \right),
\end{equation}
where $\ell_j$ $j=1,\dots , s$ are the Lagrange basis functions and so
$$ 
\sigma(\tau h)=h\sum_{j=1}^s\int_0^{\tau}\ell_j(\xi)\,d\xi \,\mathrm{dexp}_{\sigma_j}^{-1}\big( \bar{d}H_j   \righthalfcup\bar{\omega}_j\big).
$$


\vskip0.3cm
\noindent{\it Energy preservation}\newline
\vskip0.2cm
To prove energy preservation of the proposed method we consider the path $X(\xi h)=R_{x_0}\exp( \sigma (\xi h) ) \in G$ such that $X(0)=x_0$ and $X(h)=x_1$ and integrate along this path obtaining
$$H(x_1)-H(x_0)=\int_0^1\langle dH_{\exp(\sigma(\xi h))x_0}, \frac{d}{d \xi}\exp(\sigma(\xi h))x_0\rangle\, d\xi,$$
using that
$$\frac{d}{d \xi}\exp(\sigma(\xi h))x_0=R_{\exp(\sigma(\xi h))x_0}'\,\mathrm{dexp}_{\sigma(\xi h)}\left(\frac{d}{d \xi} \sigma(\xi h)\right),$$
and \eqref{eq:collpoly}, we obtain
$$H(x_1)-H(x_0)=\int_0^1\langle \mathrm{dexp}_{\sigma(\xi h)}^*R_{\exp(\sigma(\xi h))x_0}^{*}\,dH_{\exp(\sigma(\xi h))x_0},\sum_{j=1}^s\ell_j(\xi)\mathrm{dexp}_{\sigma_j}^{-1}\left( \bar{d}H_j   \righthalfcup\bar{\omega}_j \right)\rangle\, d\xi$$
and further
$$H(x_1)-H(x_0)=\sum_{j=1}^sb_j\langle\int_0^1\frac{\ell_j(\xi)}{b_j} \left( \mathrm{dexp}_{\sigma_j}^{-1} \right)^*\mathrm{dexp}_{\sigma(\xi h)}^*R_{\exp(\sigma(\xi h))x_0}^{*}\,dH_{\exp(\sigma(\xi h))x_0}\, d\xi,\bar{d}H_j   \righthalfcup\bar{\omega}_j\rangle,$$
so that finally
$$H(x_1)-H(x_0)=\sum_{j=1}^sb_j\langle\bar{d}H_j,\bar{d}H_j\righthalfcup\bar{\omega}_j \rangle =0.$$

\vskip0.3cm
\noindent{\it Order}\newline
\vskip0.2cm
Considering the change of variables $x(\xi h)=\exp(\Omega(\xi h))x_0$ $\xi\in[0, 1]$, and $x$ the solution of \eqref{eq:diffeqbivector}, by differentiation we get the following differential equation for $\Omega$:
$$\frac{d}{d \xi} \Omega=\mathrm{dexp}_{\Omega}^{-1}\left( R_{\exp(-\Omega(\xi h))}'(dH \righthalfcup \omega) \right).$$
Depending on the choice of quadrature points and weights, the collocation method \eqref{eq:collmethod1}, \eqref{eq:collmethod2} and \eqref{eq:DGho} approximates $\Omega$ at a certain order $p$, this suffices to guarantee that the overall method attains the same order, see also \cite{cohen11lep} and \cite{zanna99car}.

\begin{example}
Consider the collocation points $c_{1,2}=\frac{1}{2}\mp \frac{\sqrt{3}}{6}$ (the nodes of the Gauss quadrature) and let $\ell_1(\xi )=(\xi-c_2)/(c_1-c_2)$ and $\ell_2(\xi )=(\xi-c_1)/(c_2-c_1)$ be the corresponding Lagrange basis functions, let  
$$\sigma(\xi h)=\int_0^{\xi}\ell_1(\tau )\,d\tau \,\mathrm{dexp}_{\sigma_1}^{-1}\left( \bar{d}H_1   \righthalfcup\bar{\omega}_1\right)+\int_0^{\xi}\ell_2(\tau )\,d\tau \,\mathrm{dexp}_{\sigma_2}^{-1}\left( \bar{d}H_2   \righthalfcup\bar{\omega}_2\right),$$
and with 
$\bar{d}H_j $, $j=1,2$, given by \eqref{eq:DGho},with $b_1=b_2=1$, and $\bar{\omega}_j$, $j=1,2$, given by \eqref{eq:omegaj}.  Then the one step method 
$$x_1=\exp (\sigma( h))\,x_0$$
is a symmetric energy-preserving method of order $4$.
\end{example}

\section{Examples and numerical experiments}

In the numerical experiments we consider  the equations of the attitude rotation of a free rigid body in unit quaternions and a problem of elasticity the equations  for pseudo-rigid bodies on the cotangent bundle of ${GL}(3)$ (or $SL(3)$). 

\subsection{Attitude of a free rigid body}

The set
$$S^3=\{ \mathbbm{q}\in \mathbf{R}^4 \,\, | \,\, \| \mathbbm{q} \|^2=1 \},\quad \mathbbm{q}=[q_0,\mathbf{q}]^T,\, \mathbf{q}\in \mathbf{R}^3$$
with the quaternion product
$$\mathbbm{p}\cdot \mathbbm{q} :=[p_0q_0-\mathbf{p}^T\mathbf{q},\,p_0\mathbf{p}+q_0\mathbf{q}+\mathbf{p}\times \mathbf{q}].$$
is a Lie group. The corresponding Lie algebra is
$$\mathfrak{s}^3:=\{[0,\mathbf{v}]\in\mathbf{R}^4\,|\, \mathbf{v}\in\mathbf{R}^3\},$$
and can be identified with $\mathbf{R}^3$.

We consider the representation of the Euler attitude equations for the free rigid body in unit quaternions:
$$\dot{\mathbbm{q}}=f(\mathbbm{q})\cdot \mathbbm{q},\qquad f(\mathbbm{q})=\mathbbm{q}\cdot \mathbbm{v} \cdot \mathbbm{q}_c,$$
and 
$$\mathbbm{v}=[0\, , \mathbf{v}],\quad \mathbf{v}=\frac{1}{2}\,\mathbb{I}^{-1}\mathcal{E}(\mathbbm{q}_c)\mathbf{m}_0,$$
and $\mathbb{I}$ is the inertia tensor (a $3\times 3$ fixed, diagonal matrix).
The Euler-Rodriguez map $\mathcal{E}:S^3\mapsto SO(3)$ is defined by
$$\mathcal{E}(\mathbbm{q}):=I_3+2q_0\hat{\mathbf{q}}+2\hat{\mathbf{q}}^2,$$
where $I_3$ is the $3\times 3$ identity matrix and $\hat{\mathbf{q}}$ is defined as follows by means of the components of $\mathbf{q}$,
$$\hat{\mathbf{q}}:=\left[\begin{array}{ccc}
0 & -q_3 & q_2\\
q_3 & 0 &-q_1\\
-q_2 & q_1 & 0
\end{array}
\right].$$
The energy function preserved along $\mathbbm{q}(t)$ is
$$H(\mathbbm{q})=\frac{1}{2}\,\mathbf{m}_0^T\mathcal{E}(\mathbbm{q})\mathbb{I}^{-1}\mathcal{E}(\mathbbm{q}_c)\mathbf{m}_0.$$
We can choose the Euclidean inner product on $\mathfrak{s}^3:=\{[0,\mathbf{v}]\in\mathbf{R}^4\,|\, \mathbf{v}\in\mathbf{R}^3\}$, to play the role of the Riemannian metric on $S^3$, and 
the corresponding Riemannian gradient is 
$$\mathrm{grad}\, H=(I_4-\mathbbm{q}\mathbbm{q}^T)\,\nabla H,$$
with $I_4$ the $4\times 4$ identity.
It follows that the bivector $ \omega = \frac{\grad\; H \wedge F}{\|\grad\; H\|^2} $ can be expressed by the $4\times 4$ rank-$2$ matrix
 $$    
     \omega_R(\mathbbm{q})=\frac{\mathbf{\xi}\mathbf{\gamma}^T-\mathbf{\gamma}\mathbf{\xi}^T}{\|\mathbf{\gamma}\|^2},\,\, \mathbf{\xi},\mathbf{\gamma}\in \mathfrak{s}^3,
$$
where
$\mathbf{\xi}=f(\mathbbm{q})$ and $\mathbf{\gamma}=\left.\mathrm{grad} H\right|_{\mathbbm{q}}\cdot \mathbbm{q}_c $. As a discrete bivector to construct the method we choose  
$$\bar{\omega}(\mathbbm{q},\mathbbm{q}')=\omega_R(\bar{\mathbbm{q}}),\quad \bar{\mathbbm{q}}=\exp(\eta/2)\mathbbm{q},\quad\eta=\log(\mathbbm{q}'\cdot\mathbbm{q}_c).$$

Identifying $\mathfrak{s}^3$ with its dual, and using the Gonzales trivialised discrete differential we get
$$\overline{\mathrm{grad}\, H}(\mathbbm{q},\mathbbm{q}')=
\mathbf{\gamma} 
+\frac{H(\mathbbm{q}')-H(\mathbbm{q})-
  \mathbf{\gamma}^T\eta 
}{\|\eta\|^2}\eta.$$

The energy-preserving symmetric Lie group method is
$$\mathbbm{q}'=\exp\big(h\,\bar{\omega}(\mathbbm{q},\mathbbm{q}')\overline{\mathrm{grad}\, H}(\mathbbm{q},\mathbbm{q}')\big)\cdot \mathbbm{q}.$$

Alternatively we can use the averaged trivialised discrete differential obtained averaging the Riemannian gradient as
\begin{equation}
\label{avfTDD}
\overline{\mathrm{grad}\, H}(\mathbbm{q},\mathbbm{q}')=
\int_0^1 \mathbf{\gamma}(\xi)\, d\xi, \quad \mathbf{\gamma}({\xi})=\left.\mathrm{grad} H\right|_{\mathbbm{q}(\xi)}\cdot \mathbbm{q}_c(\xi), 
\end{equation}
and $\mathbbm{q}(\xi)=\exp(\xi\, \log(\mathbbm{q}'\cdot \mathbbm{q}_c))\cdot \mathbbm{q}$.

\begin{figure}[htp]\centering
\begin{tabular}{cc}
\includegraphics[width=50mm]{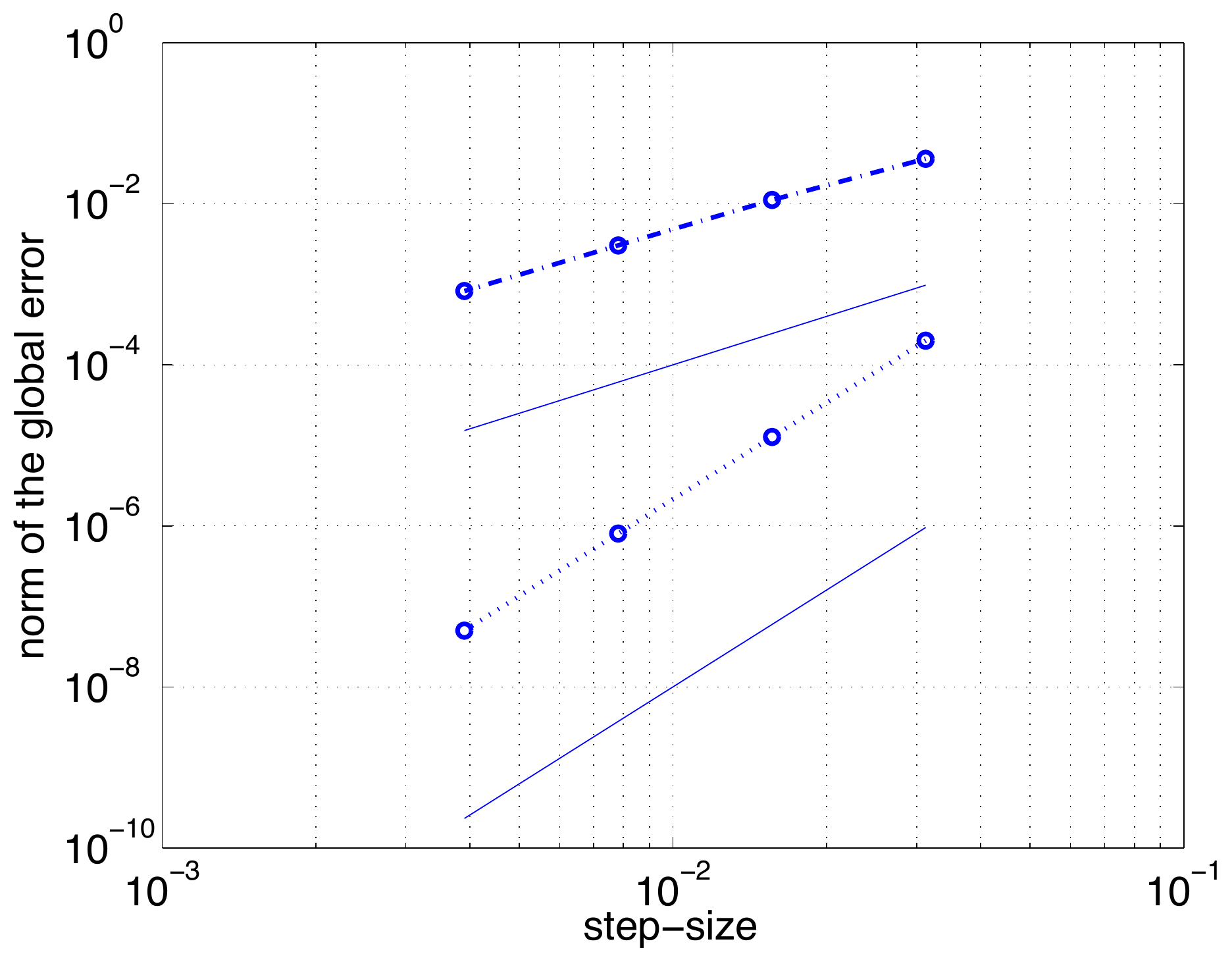} 
\includegraphics[width=50mm]{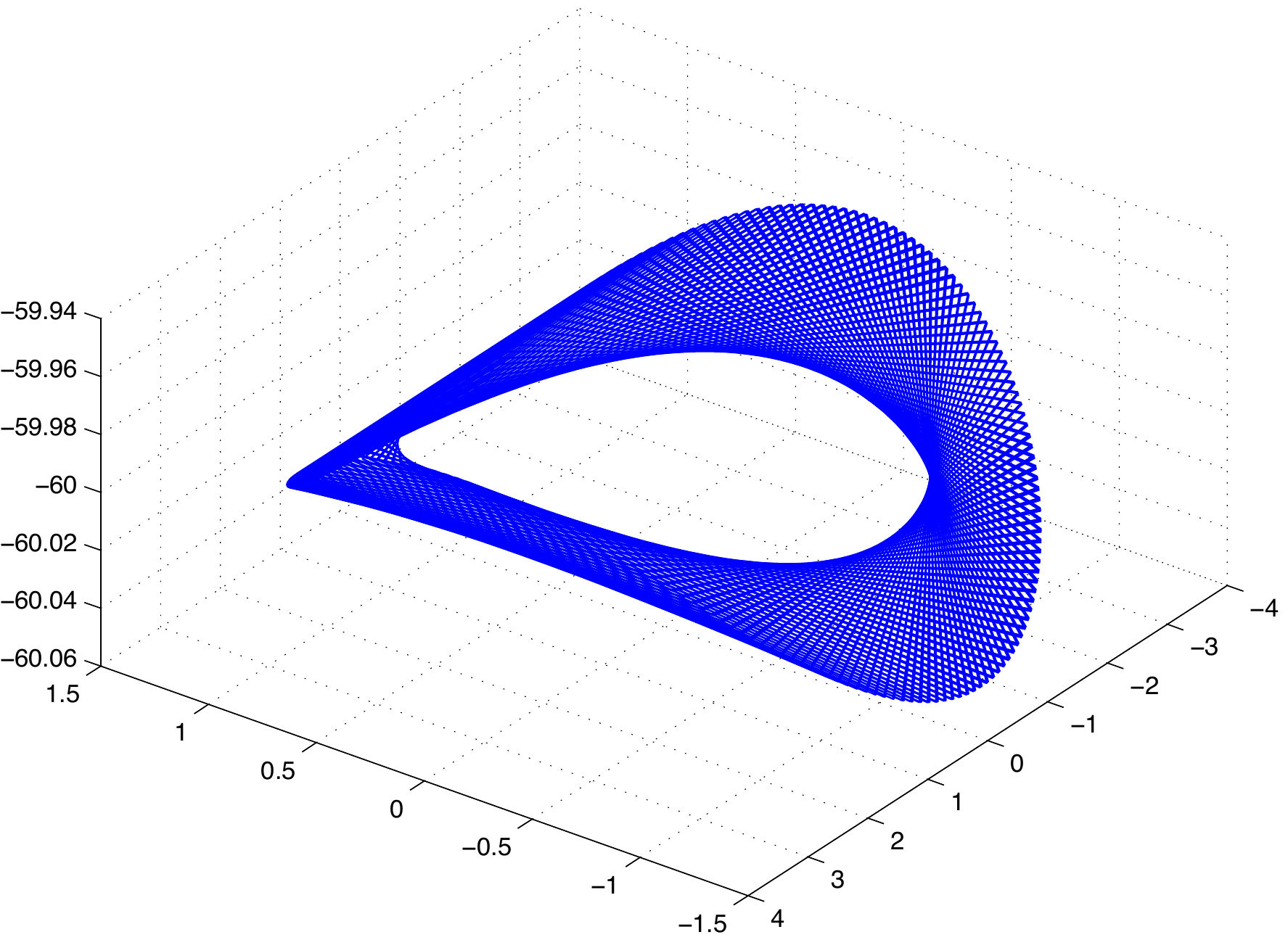}\\
\end{tabular}
\caption{Rigid body problem, $\mathbbm{q}_0= [1,0,0,0]^T$. Inertia tensor $\mathbbm{I}=\mathrm{diag}(1,5,60)$, $\mathbf{m}_0=\mathbbm{I}\mathbf{v}_0$ and $\mathbf{v}_0=[1, 0.5,-1]^T$. (Left) Order of the energy preserving methods, step-size versus norm of the global error: dashed-dotted line second order symmetric energy-preserving method based on the TDD \eqref{avfTDD}; dotted line symmetric energy-preserving method of order $4$ (section~\ref{sec:4}); the solid lines are reference lines of order $2$ and $4$.  (Right) momentum vector for the energy-preserving symmetric Lie group method of order $4$, time interval $[0,50]$, step-size $h=2^{-4}$. \label{fig:EPhigh}}
\end{figure}

In figure~\ref{fig:EPhigh} we apply a symmetric energy-preserving method of order 2 and 4 to the free rigid body problem in the formulation presented in this section, and we compare them with an explicit Lie group method of order 2 based on Heun's Runge-Kutta formula. We have used the discrete trivialized differential (\ref{avfTDD}) in the Lie group method of order two and the formulation outlined in section \ref{sec:4} for the $4$-th order variant of the method. The collocation points for this method are $c_{1,2}=\frac{1}{2}\mp \frac{\sqrt{3}}{6}.$ The integrals are approximated with accurate quadrature formulae. 

\subsection{Pseudo-rigid bodies}

We consider the Hamiltonian equations describing rotating homogeneous elastic rigid bodies \cite{lewis94caf}, \cite{lewis90nso}. A pseudo-rigid body is a three dimensional elastic body whose deformation gradient $F=F(t)\in GL^{+}(3)$ is assumed to be constant throughout the body $\mathcal{B}$: for all $X\in \mathcal {B}\subset \mathbf{R}^3$,  the deformation is always given by  the matrix vector product 
$F\,X$ and $F$ is not depending on $X$. In the incompressible case $F\in SL(3)$. The configuration space of a pseudo-rigid body is the group $G$ equal to $GL^{+}(3)$ or $SL(3)$, and the flow of the corresponding Hamiltonian equations evolves on $T^*G\approx G \times \mathfrak{g}^*$.  With the semidirect-product group structure induced by the group multiplication in $G$, $G \times \mathfrak{g}^*$ is also a Lie group.  
In our particular example  $GL^{+}(3)\times \mathfrak{gl}(3)^*$, we identify $\mathfrak{gl}(3)^*$ with its dual, 
and we use  coordinates $(F,P)$ where $F\in GL^{+}(3)$ and $P\in \mathfrak{gl}(3)$ respectively.  The corresponding Lie algebra is $\mathfrak{gl}(3)\times \mathfrak{gl}(3)^*$. We denote with $W(C)$ the stored energy function depending on the Cauchy-Green tensor $C=F^TF$, and  with $E$ the inertia tensor. The Hamiltonian function is
$$H(F,P):=K(P)+W(F^TF),\qquad K(P)=\frac{1}{2}\, \langle P, PE^{-1} \rangle,$$
where $\langle P ,V\rangle=\mathrm{tr}(P^TV)$ is the standard matrix inner product, i.e. the duality pairing between between $TG$ and $T^*G$. The canonical Hamilton's equations take the form
\begin{eqnarray}
\label{eq:prb}
\dot{F}&= & \frac{\delta H}{\delta P},\\
\dot{P}&= & -\frac{\delta H}{\delta F},
\end{eqnarray}
where 
$(\frac{\delta H}{\delta F},\frac{\delta H}{\delta P})=(-2F\, \nabla W(F^TF), PE^{-1})$.

In our experiments we consider a St Venant-Kirchhoff material leading to a stored energy function
$$W(C)=\frac{1}{2}\, \lambda (\mathrm{tr}(C-I))^2+\mu\, \mathrm{tr}((C-I)^2),$$
with $\mathrm{tr}$ denoting the trace operator and $\lambda$ and $\mu$ the Lam{\'e} constants (such that $\mu >0$ and $3\lambda+2\mu >0$, $\mu=1$ and $\lambda=\frac{1}{3}$ in the experiments).  

The bivector $\omega$ is represented in this case simply by the $6\times 6$ inverse Darboux matrix
$$J^{-1}=\left[
\begin{array}{cc}
O & I\\
-I & O
\end{array}
\right].$$

The semi-direct product group multiplication in $G\times \mathfrak{g}^*$ is 
$$(F_1,P_1)\cdot (F_2,P_2)=(F_1F_2, P_1+\mathrm{Ad}^{*}_{F_1^{-1}}\,P_2),$$
and as matrix operation $\mathrm{Ad}^{*}_{F_1^{-1}}\,P_2=F_1^{-T}P_2F_1^{-T}$,
we denote with $\mathrm{Exp}$ and $\mathrm{Log}$ the exponential and logarithm between $G\times \mathfrak{g}^*$ and $\mathfrak{g}\times \mathfrak{g}^*$ respectively, these are defined by: 
$$\mathrm{Exp}((\xi, \mu))=(\exp(\xi),\mathrm{dexp}_\xi(\mu)),\qquad \mathrm{Log}((g,\sigma))=(\log(g),(\mathrm{dexp}^*_{-\log(g)})^{-1}(\sigma)),
$$
where $\exp$ and $\log$ are the corresponding maps for $G$.
We consider
$$\eta=(\eta_F,\eta_P)=\mathrm{Log}((F_1,P_1)\cdot (F_0,P_0)^{-1}),
$$
and denote with $(\gamma_1,\gamma_2):=R^*_{(\bar{F},\bar{P})}\, \left. dH\right|_{(\bar{F},\bar{P})}\in \mathfrak{gl}(3)^*\times \mathfrak{gl}(3)$ the trivialized differential in the point $(\bar{F},\bar{P})=\mathrm{Exp}(\eta_F/2,\eta_P)(F_0,P_0)$, we have
$$\gamma_1=\left.\frac{\delta H}{\delta F}\right|_{(\bar{F},\bar{P})}-\mathrm{ad}_{\gamma_2}^*(\bar{P}), 
\qquad\gamma_2=\left. \frac{\delta H}{\delta P}\right|_{(\bar{F},\bar{P})}\bar{F}^{-1},$$
and in matrix form $\mathrm{ad}_{\gamma_2}^*(\bar{P})=\gamma_2^T\bar{P}-\bar{P}\gamma_2^T$.
The trivialised discrete differential \eqref{eq:tddgmp} becomes in coordinates
$$
\bar{d}H_{((F,P),(F',P'))}=(\gamma_1,\gamma_2)+\alpha (\eta_F,\eta_P),\qquad \alpha:=\frac{H(F',P')-H(F,P)-\langle \gamma_1,\eta_F\rangle-\langle \eta_P,\gamma_2\rangle}{\|\eta\|^2},
$$
and
where the metric is deduced by the standard matrix inner product,
$$\|\eta\|^2:=\mathrm{tr}(\eta_F^T\eta_F)+\mathrm{tr}(\eta_P^T\eta_P).$$
The energy preserving method can then be formulated as 
\begin{equation}
\label{eq:syEPLGM}
(F_1,P_1)=\mathrm{Exp}(h(\gamma_2+\alpha \eta_P,-\gamma_1-\alpha \eta_F)) (F_0,P_0).
\end{equation}

\begin{figure}[htp]\centering
\begin{tabular}{cc}
\includegraphics[width=50mm]{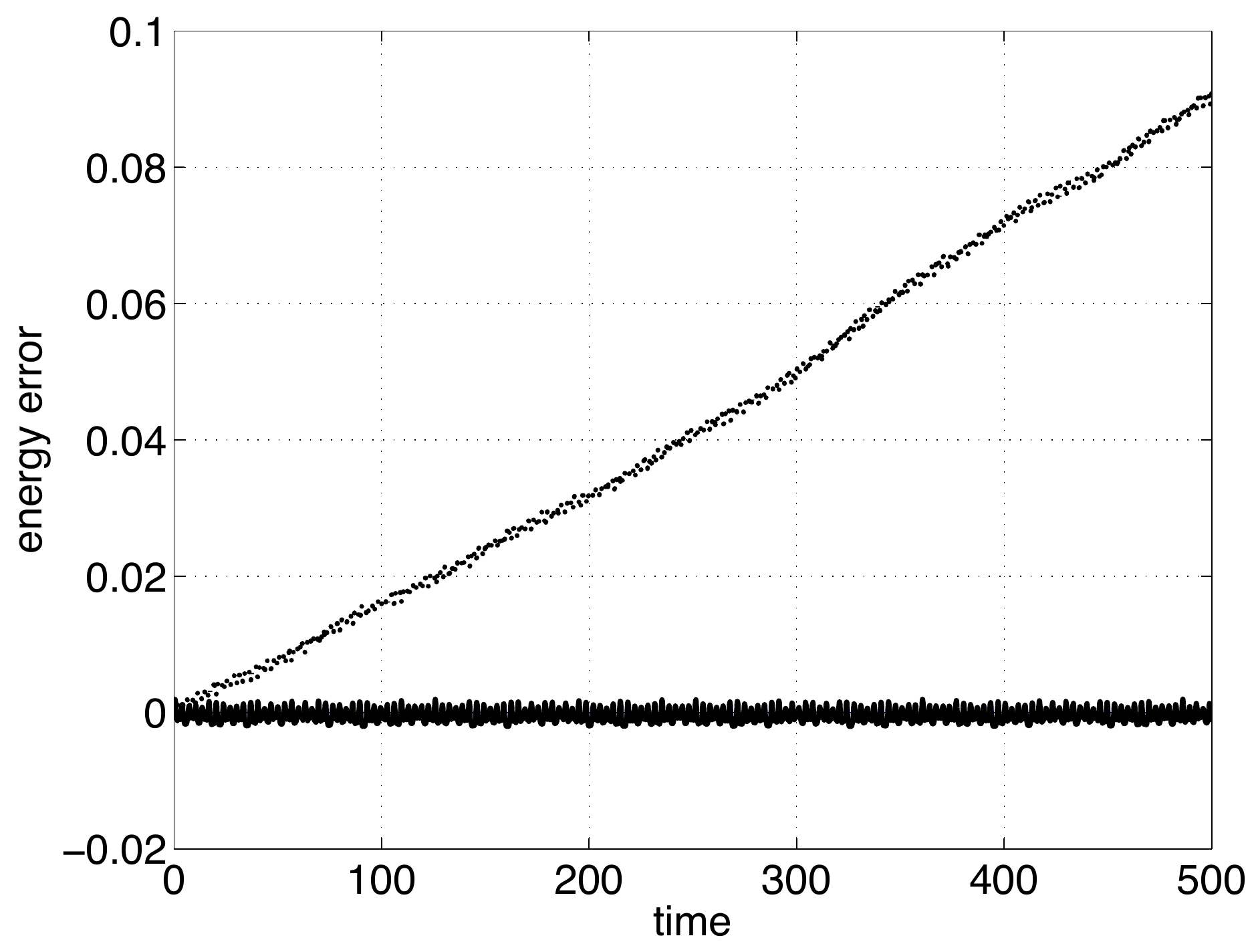} &
\includegraphics[width=50mm]{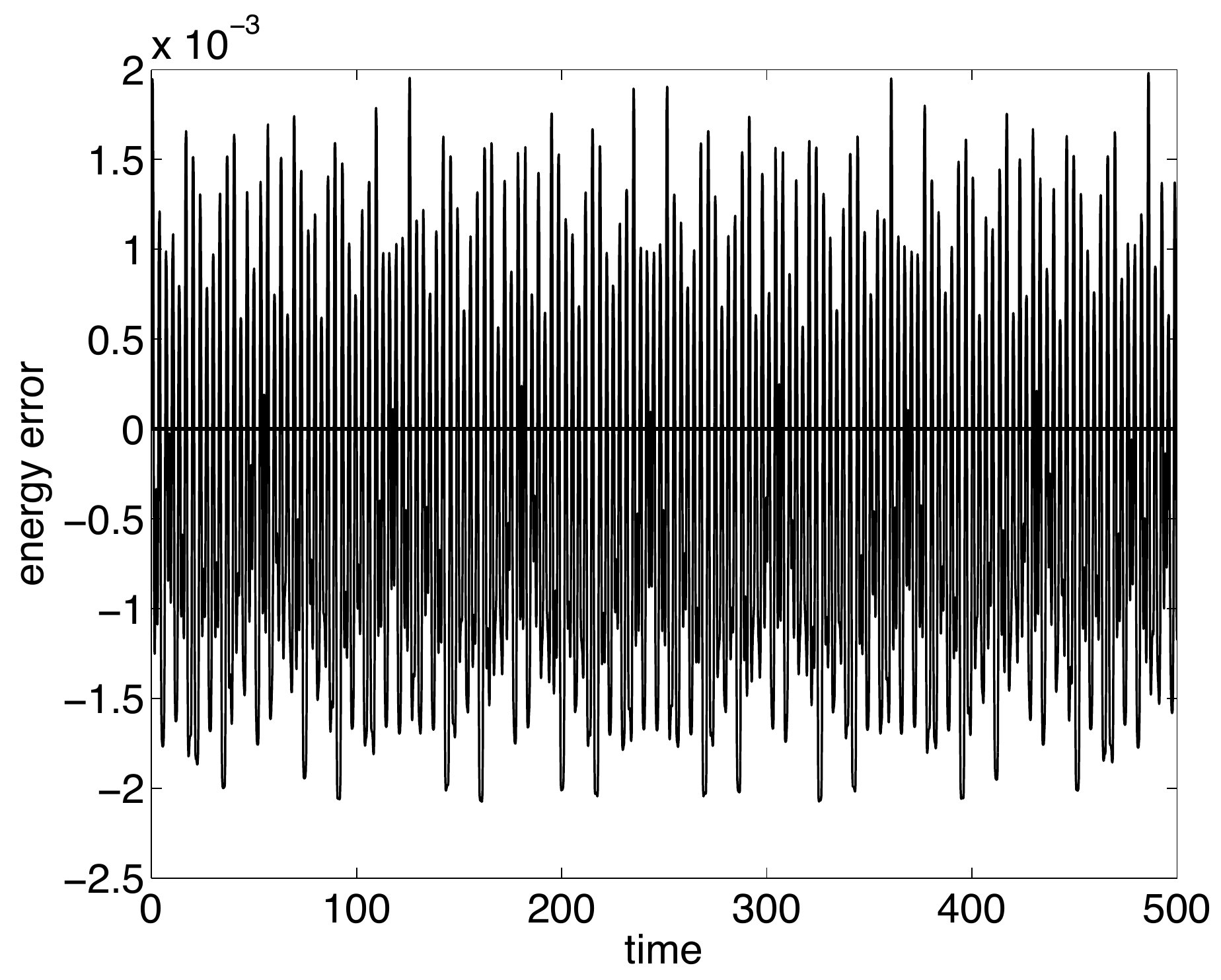}\\
\includegraphics[width=50mm]{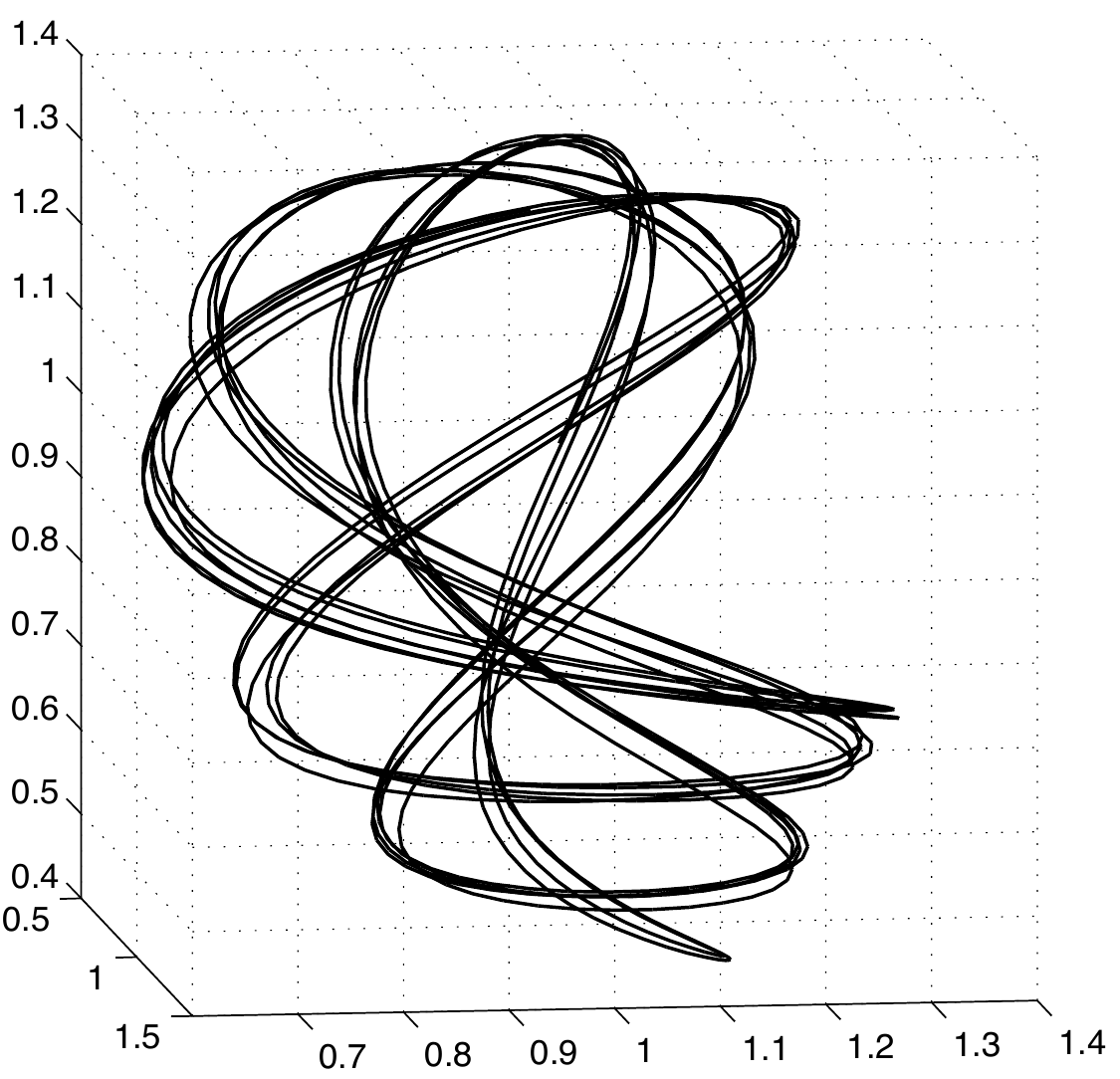} &
\includegraphics[width=50mm]{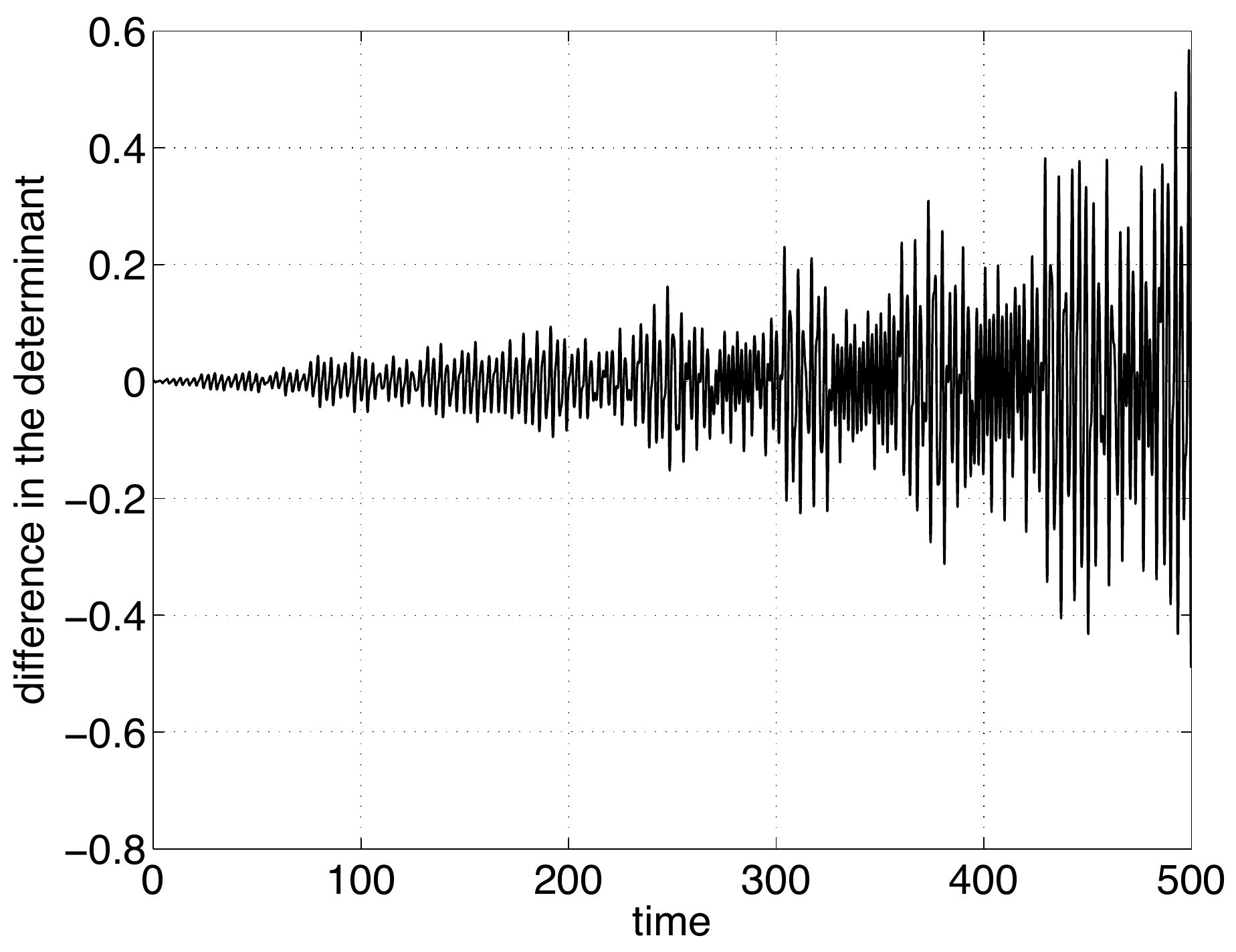}\\
\end{tabular}
\caption{Pseudo-rigid body problem, $F_0=I$, $P_0=\mathrm{diag}(0.2575, 0.8407, 0.2543)$, 
integration interval $[0,500]$. Lam{\'e} parameters $\lambda=\frac{1}{3}$, $\mu=1$. $E=\mathrm{diag}(1,2,3)$. (Top left) energy error explicit Lie group method (Heuns method), dotted line, $h=1/32$, versus symmetric order 2 Lie group method, solid line, $h=1/16$. (Top right) energy error symmetric order 2 Lie group method (energy error $10^{-3}$), versus the energy preserving Lie group method (energy error $10^{-12}$, straight line), $h=1/16$ for both methods. (Bottom left) deformation of the vector $[1,1,1]^T$ under the transformation $F(t)$, first $1000$ steps, $h=1/16$, energy-preserving method (with a different $P_0$). (Bottom right) difference of the determinants of $F$ for the symmetric Lie group method and the symmetric, energy-preserving Lie group method: $\mathrm{det}(F_n^{sym})-\mathrm{det}(F_n^{EP})$, $n=1,\dots ,8000$.  \label{fig:EP}}
\end{figure}

In figure~\ref{fig:EP} we report the results of a simulation for this problem. We compare an explicit Lie group method (Heun's method), a symmetric Lie group method and the symmetric energy-preserving Lie group method presented in this section. The symmetric Lie group method is obtained by setting $\alpha=0$ in \eqref{eq:syEPLGM}.  All methods are Lie group methods of order $2$. We have halved the step-size for the explicit second order Lie group method, to avoid instability. All Lie group methods have the property that $\det (F_n)>0$ for all $n$.  The explicit Lie group method  fails to preserve the energy of the problem (figure~\ref{fig:EP} top left), the symmetric Lie group methods have both a much smaller energy error (figure~\ref{fig:EP} top right), the symmetric energy preserving Lie group method preserves the energy to very high precision. We report here experiments with diagonal initial values for $F$ and $P$. We have performed experiments also with non diagonal initial values obtaining similar results. The performance of the methods is relying on the accurate computation of matrix functions, and in particular matrix logarithms. In figure~\ref{fig:EP} (bottom right) we show the difference in the determinants of $F$ for the two symmetric methods: the symmetric one denoted (sym) and the symmetric and energy preserving one denoted (EP). We plot $\mathrm{det}(F_n^{sym})-\mathrm{det}(F_n^{EP})$ for $n=1,\dots ,8000$, this measures how the two numerical solutions depart from each other with time.


\section*{Acknowledgments.}
This research was supported by a Marie Curie International Research Staff Exchange
Scheme Fellowship within the 7th European Community Framework Programme.
The authors would like to acknowledge the support from the GeNuIn
Applications and SpadeAce projects funded by the Research Council of Norway, and
most of the ideas arise while the authors were visiting Massey University,
Palmerston North, New Zealand and La Trobe University, Melbourne, Australia.

 \bibliographystyle{plain}
 \bibliography{geom_int,mybib,mybibE}

\begin{thebibliography}{10}

\bibitem{adler02nmo}
R.~L. Adler, J.~P. Dedieu, J.~Y. Margulies, M.~Martens, and M.~Shub.
\newblock Newton's method on {R}iemannian manifolds and a geometric model for
  the human spine.
\newblock {\em IMA Journal of Numerical Analysis}, 22(3):359--390, 2002.

\bibitem{brugnano10hbv}
Luigi Brugnano, Felice Iavernaro, and Donato Trigiante.
\newblock Hamiltonian boundary value methods (energy preserving discrete line
  integral methods).
\newblock {\em JNAIAM. J. Numer. Anal. Ind. Appl. Math.}, 5(1-2):17--37, 2010.

\bibitem{celledoni12per}
E.~Celledoni, V.~Grimm, R.I. McLachlan, D.I. McLaren, D.~O'Neale, B.~Owren, and
  G.R.W. Quispel.
\newblock Preserving energy resp. dissipation in numerical pdes using the
  "average vector field" method.
\newblock {\em Journal of Computational Physics}, 231(20):6770--6789, 2012.
\newblock cited By (since 1996) 0.

\bibitem{celledoni12ait}
E.~Celledoni, H.~Marthinsen, and B.~Owren.
\newblock An introduction to {L}ie group integrators -- basics, new
  developments and applications.
\newblock {\em Journal of Computational Physics}, 2013.
\newblock To appear.

\bibitem{christiansen11tis}
S.H. Christiansen, H.Z. Munthe-Kaas, and B.~Owren.
\newblock Topics in structure-preserving discretization.
\newblock {\em Acta Numerica}, 20:1--119, 2011.

\bibitem{cohen11lep}
David Cohen and Ernst Hairer.
\newblock Linear energy-preserving integrators for {P}oisson systems.
\newblock {\em BIT}, 51(1):91--101, 2011.

\bibitem{dahlby11agf}
M.~Dahlby and B.~Owren.
\newblock A general framework for deriving integral preserving numerical
  methods for {PDE}s.
\newblock {\em SIAM J. Sci. Comput.}, 33(5):2318--2340, 2011.

\bibitem{gonzalez96tia}
O.~Gonzalez.
\newblock Time integration and discrete {H}amiltonian systems.
\newblock {\em J. Nonlinear Sci.}, 6:449--467, 1996.

\bibitem{hairer10epc}
E.~Hairer.
\newblock Energy-preserving variant of collocation methods.
\newblock {\em Journal of Numerical Analysis, Industrial and Applied
  Mathematics}, 5(1-2):73--84, 2010.

\bibitem{iserles00lgm}
A.~Iserles, H.~Z. Munthe-Kaas, S.~P. N{\o}rsett, and A.~Zanna.
\newblock Lie-group methods.
\newblock {\em Acta Numerica}, 9:215--365, 2000.

\bibitem{lewis90nso}
D.~Lewis and J.~C. Simo.
\newblock Nonlinear stability of rotating pseudo-rigid bodies.
\newblock {\em Proc. R. Soc. Lond. A}, 427(1873):281--319, 1990.

\bibitem{lewis94caf}
D.~Lewis and J.~C. Simo.
\newblock Conserving algorithms for the dynamics of {H}amiltonian systems of
  {L}ie groups.
\newblock {\em J. Nonlinear Sci.}, 4:253--299, 1994.

\bibitem{mclachlan99giu}
R.~I. McLachlan, G.~R.~W. Quispel, and N.~Robidoux.
\newblock Geometric integration using discrete gradients.
\newblock {\em Phil. Trans. Royal Soc. A}, 357:1021--1046, 1999.

\bibitem{zanna99car}
A.~Zanna.
\newblock Collocation and relaxed collocation for the {F}er and the {M}agnus
  expansions.
\newblock {\em SIAM J. Numer. Anal.}, 36(4):1145--1182, 1999.

\end{thebibliography}

\end{document}